\documentclass{mcom-l}

\newcommand{\Z}{\mathbb{Z}}
\newcommand{\N}{\mathbb{N}}

\newcommand{\PP}{\mathbb{P}}

\DeclareMathOperator{\ord}{ord}
\DeclareMathOperator{\lcm}{lcm}
\DeclareMathOperator{\lev}{lev}
\DeclareMathOperator{\integ}{int}

\usepackage{amssymb}
\usepackage{tikz}
\usetikzlibrary{arrows,backgrounds}

\copyrightinfo{}{}

\newtheorem{theorem}{Theorem}[section]
\newtheorem{lemma}[theorem]{Lemma}
\newtheorem{corollary}[theorem]{Corollary}

\theoremstyle{definition}
\newtheorem{definition}[theorem]{Definition}
\newtheorem{example}[theorem]{Example}
\newtheorem{problem}[theorem]{Problem}

\theoremstyle{remark}
\newtheorem{remark}[theorem]{Remark}

\numberwithin{equation}{section}

\begin{document}

\title[Factorization and Tetration]{A reduction of integer factorization to modular tetration}


\author[M. Hittmeir]{Markus Hittmeir}
\address{Hellbrunnerstra{\ss}e 34, A-5020 Salzburg}
\curraddr{}
\email{markus.hittmeir@sbg.ac.at}
\thanks{The author is supported by the Austrian Science Fund (FWF): Project F5504-N26,
which is a part of the Special Research Program "Quasi-Monte Carlo Methods: Theory
and Applications".}


\date{}

\dedicatory{}

\begin{abstract}
Let $a,k\in\N$. For the $k-1$-th iterate of the exponential function $x\mapsto a^x$, also known as tetration, we write
\[
^k a:=a^{a^{.^{.^{.^{a}}}}}.
\]
In this paper, we show how an efficient algorithm for tetration modulo natural numbers $N$ may be used to compute the prime factorization of $N$. In particular, we prove that the problem of computing the squarefree part of integers is deterministically polynomial-time reducible to modular tetration.
\end{abstract}

\maketitle

\section{Introduction}
The fastest currently known algorithms for computing the prime factorization of natural numbers achieve sub-exponential runtime complexities. In practice, methods such as the Number Field Sieve, Pomerance's Quadratic Sieve and Lenstra's Elliptic Curve Method are used for factoring large integers. The monographs \cite{Rie} and \cite{Wag} provide an overview on the variety of methods. Amongst these approaches, Shor's algorithm \cite{Sho} takes in an exceptional position. It runs in polynomial-time, but relies on quantum computing. Since powerful quantum computers do not yet exist, Shor's result has remained theoretical. Moreover, all methods mentioned above are either conditional in some sense or rely on probabilistic and heuristic arguments. The lowest rigorously established complexity bounds are still exponential in the length of the input number. In 1977, Strassen \cite{Str} published a deterministic factorization method running in about $\widetilde{O}(N^{1/4})$ computational steps. Since then, there have been improvements due to Bostan, Gaudry and Schost \cite{BosGauSch}, Costa and Harvey \cite{CosHar} and myself \cite{Hit}.  

This paper concerns the relations of the integer factorization problem ($IFP$) to other computational tasks. Factoring natural numbers is widely assumed to be hard and very little is known about the possibility of non-quantum polynomial-time solutions. On the one hand, a fast factorization method would allow to efficiently solve a variety of other problems, such as deciding quadratic residuosity or computing Carmichael's function. On the other hand, only a few problems are known to which $IFP$ is polynomial-time reducible. One example is the discrete logarithm problem. An efficient technique for solving it would yield an efficient algorithm for factoring integers. Up to now, it is unclear whether the converse statement is true or not. The evaluation of Carmichael's function is a second example and, therefore, actually polynomial-time equivalent to factorization.

Our main goal is to show that $IFP$ is related to yet another algorithmic problem. Let $a,k\in\N$. The first four hyperoperations are as follows:
\begin{enumerate}
\item{Addition: $a+k=a+\underbrace{1+1+\cdots+1}_k$ may be considered as the $k$-th iterate of the successor function evaluated at $a$.}
\item{Multiplication: $a\cdot k=\underbrace{a+a+\cdots +a}_k$ may be considered as the $k-1$-th iterate of $x\mapsto a+x$ evaluated at $a$.}
\item{Exponentiation: $a^k=\underbrace{a\cdot a\cdots a}_k$ may be considered as the $k-1$-th iterate of $x\mapsto a\cdot x$ evaluated at $a$.}
\item{Tetration: $^k a=\underbrace{a^{a^{.^{.^{.^{a}}}}}}_k$ may be considered as the $k-1$-th iterate of $x\mapsto a^x$ evaluated at $a$.}
\end{enumerate}
It is clear how this list could be continued ad infinitum. 

We are interested in the complexity of performing hyperoperations modulo natural numbers $N$. Let $a,k\in\Z/N\Z$. Computing $a+k\pmod N$ takes $O(\log N)$ bit operations. The bit-complexity for computing $a\cdot k\pmod N$ is in $O(\textsf{M}_{\integ}(\log N))$, where $\textsf{M}_{\integ}(j)$ denotes the cost for multiplying two $\lceil j\rceil$-bit integers. Due to \cite{HarHoeLec}, $\textsf{M}_{\integ}(j)$ can be bounded by 
$
O(k\log k\cdot 8^{\log^*k}),
$
where $\log^*j$ stands for the iterated logarithm. Now let $a\in\Z/N\Z$ and $k\in\N$. Then the value $a^k \pmod N$ may be computed in $O(\textsf{M}_{\integ}(\log N)\cdot\log k)$ bit operations by using the Square-and-Multiply procedure. We see that the runtime complexities for the first three hyperoperations are unconditional in the sense that they do not use any additional information about the modulus $N$ and, in particular, do not require knowledge about its prime factorization. In addition to that, the algorithms run in polynomial-time with respect to the input triple $(a,k,N)$. Naturally, these observations give rise to the problem of constructing such an algorithm for computing the value of $^k a\pmod N$. In the following, we will call this the modular tetration problem ($MTP$).

Our contribution is to present a reduction from $IFP$ to $MTP$ and demonstrate its efficiency. In Sections 2 and 3, we will discuss the work established in \cite{BlaBor} and show that it directly implies the existence of a deterministic polynomial-time reduction from $MTP$ to $IFP$. In Section 4, we will consider our main idea for the reverse statement. In order to investigate the efficiency of the reduction, we introduce the notion of the \emph{level} modulo natural numbers and prove some of its basic properties in Section 5. Two results arise from this study:
\begin{itemize}
\item{Theorem \ref{t3} in Section 6 is the main result of this paper. It states that the problem of computing the squarefree part of integers is deterministically polynomial-time reducible to modular tetration. Integer squarefree decomposition is important due to its relation to one of the main problems of computational algebra, namely determining the ring of integers of an algebraic number field. In Corollary \ref{c5}, we obtain that this task is also deterministically polynomial-time reducible to $MTP$. }
\item{Section 7 contains an analysis of the general reduction from $IFP$ to $MTP$. Based on Theorem \ref{t4}, we provide heuristic arguments for the effective applicability to most integers and discuss the nature of exceptions.}
\end{itemize}
In Theorem \ref{t2}, we furthermore prove that $MTP$ is deterministically polynomial-time reducible to the discrete logarithm problem.

\section{Preliminaries}
In the following section, we will recall elementary and already known results concerning tetration modulo natural numbers. The most important proofs are carried out completely; however, the interested reader may find a more detailed exposition of all of the subsequent notions and ideas in \cite{BlaBor}. Throughout the paper, let $\PP$ denote the set of prime numbers, let $\N_0=\N\cup\{0\}$ and set $\Z_n^*=(\Z/n\Z)^*$ for $n\in\N$. Furthermore, put $^0 a=1$ to complete the definition of $^k a$ for every $k\in\N_0$.

{\definition{Let $n\in\N$. The smallest number $m$ such that 
\[
a^m\equiv 1 \mod n
\]
for every $a\in\Z_n^*$ is called the \emph{group exponent of $\Z_n^*$}. The mapping $\lambda:\N\rightarrow \N$ assigning every natural number $n$ to the group exponent of $\Z_n^*$ is called \emph{Carmichael function}.}}

{\remark{Let $n,d\in\N$. It is easy to check that Carmichael's function satisfies the following properties:
\begin{itemize}
\item{$\lambda(n)<n$ for $n\geq 2$,}
\item{$\lambda(n)$ is even for $n\geq 3$,}
\item{$\lambda(n)\leq n/2$ for even $n$,}
\item{$\lambda(d)\mid\lambda(n)$ for $d\mid n$.}
\end{itemize}
}\label{r1}}

{\lemma{The following holds:
\begin{itemize}
\item{$\lambda(1)=\lambda(2)=1$, $\lambda(4)=2$,}
\item{$\lambda(2^k)=2^{k-2}$ for $k>2$,}
\item{$\lambda(p^k)=p^{k-1}(p-1)$ for $p\geq 3$ prime and $k\in\N$,}
\item{$\lambda(p_1^{\alpha_1}\cdots p_r^{\alpha_r})=\lcm\{\lambda(p_1^{\alpha_1}),...,\lambda(p_r^{\alpha_1})\}$ for $p_i\in\PP$ and $k_i\in\N$, $i=1,...,r$.}
\end{itemize}
}\label{l1}}

\begin{proof}
These properties of Carmichael's function are well-known. A proof can be found in \cite[pp. 53-54]{Lev}.
\end{proof}

The iterates of Carmichael's function play an important role in our following considerations.

{\definition{Let $n\in\N$. We define:
\begin{itemize}
\item{$\lambda^{(0)}(n):=n$,}
\item{$\lambda^{(k)}(n):=\lambda(\lambda^{(k-1)}(n))$, for all $k\in\N$,}
\item{$H(n):=\min\{\alpha:\lambda^{(\alpha)}(n)=1\}$,}
\item{$L(0,n):=1$,}
\item{$L(t,n):=\lcm\{\lambda^{(H(n))}(n),\lambda^{(H(n)-1)}(n),...,\lambda^{(H(n)-t)}(n)\}$ for all $t\in\N$,}
\item{$L(n):=L(H(n),n)=\lcm\{n,\lambda(n),\lambda^{(2)}(n),...,\lambda^{(H(n))}(n)\}$,}
\item{$E(n):=\max\{\nu_p(n): p\mid n\}$,}
where $\nu_p(n)$ is the exponent of $p$ in the prime factorization of $n$.
\end{itemize}
}}

{\lemma{For $n\in\N$, the following holds:
\begin{enumerate}
\item{$H(n)\leq \lceil \log_2(n)\rceil$,}
\item{$E(n)\leq 2H(n)-1$ for $n>1$.}
\end{enumerate}
\label{l2}}}

\begin{proof}
A proof can be found in \cite[pp. 569-570]{BlaBor}. We carry out the argument for (1) and start by showing $H(n)\leq \log_2(n)$ for even $n$. Clearly, $H(2)=1=\log_2(2)$. For even $n\geq 4$, Remark \ref{r1} implies that $\lambda(n)$ is even and not larger than $n/2$. By induction, we derive
\[
H(n)=1+H(\lambda(n))\leq 1+\log_2(\lambda(n))\leq 1+\log_2(n/2)=\log_2(n).
\]
We now consider odd numbers $n$. We have $H(1)=0=\log_2(1)$, and for odd $n\geq 3$ we know that $\lambda(n)$ is even and strictly smaller than $n$. Applying the already shown statement for even numbers, we deduce
\[
H(n)=1+H(\lambda(n))\leq 1+\log_2(\lambda(n))<1+\log_2(n)
\]
and, since $\log_2(n)$ cannot be an integer, the claim follows.
\end{proof}

\definition{Let $n,a\in\N$. Let $V=V(a,n)$ be defined as the largest factor of $n$ coprime to $a$, and let $W=W(a,n):=n/V$. We call $n=VW$ the \emph{orthogonal decomposition of $n$ with respect to $a$}.\label{d1}}

{\lemma{Let $n,a,k,j\in\N$. If $k\geq j\geq E(n)$ and $k\equiv j \mod \lambda(n)$, then
\[
a^k\equiv a^j \mod n.
\]}\label{l3}}

\begin{proof}
This is Lemma 2.3 in \cite{BlaBor}. In the proof, the orthogonal decomposition of $n$ with respect to a is considered. It is easy to check that the congruence holds modulo $V(a,n)$ and modulo $W(a,n)$, and an application of the Chinese Remainder Theorem concludes the proof.
\end{proof}

{\lemma{Let $n,a,b\in\N$. If $a\equiv b \mod L(n)$ and $a\geq b\geq E(n)$, then
\[
^k a\,\equiv \, ^k b \mod L(n)
\]
for every $k\in\N$.
}}

\begin{proof}
This is Lemma 4.1 in \cite{BlaBor}. The proof is straight-forward and uses Lemma \ref{l3} and induction on $k$.
\end{proof}

{\theorem{Let $n,a,k,j\in\N$. If $k\geq j\geq H(n)+1$, then
\[
^ka\, \equiv\, ^ja \mod L(n).
\]}\label{t1}}

\begin{proof}
This is Lemma 3.1 in \cite{BlaBor}. For the proof, we may assume that $a,n\geq 2$. Furthermore, let $h=H(n)$. We will show by induction on $t$ that the congruence
\[
^{u+t}a\, \equiv \, ^{y+t} a \mod L(t,n)
\]
holds for all $u,y\in\N$ and all nonnegative integers $t\leq h$. For $t=0$, we have $L(0,n)=1$ and the claim follows. Now assume that, for some $t<h$, the congruence above holds for all $u,y$. This implies
\begin{align}
^{(u+1)+t}a\, \equiv \, ^{(y+1)+t} a \mod L(t,n).
\end{align}
Note that since $\lambda^{(h-t)}(n)$ divides $L(t,n)$, we also get $^{u+t} a\, \equiv \, ^{y+t} a \mod \lambda^{(h-t)}(n)$. Using (2) of Lemma \ref{l2}, we derive
\[
^{u+t} a\geq 2^{t+1}\geq 2(t+1)=2H(\lambda^{(h-t-1)}(n))>E(\lambda^{(h-t-1)}(n)).
\]
By similar arguments, one furthermore obtains $^{y+t} a>E(\lambda^{(h-t-1)}(n))$. We may apply Lemma \ref{l3} to deduce
\begin{align}
^{u+t+1}a\,=\,a^{^{u+t}a}\,\equiv\, a^{^{y+t}a}\,=\,^{y+t+1}a\mod \lambda^{(h-t-1)}(n).
\end{align}
Since $L(t+1,n)=\lcm\{L(t,n),\lambda^{h-t-1}(n)\}$, we now combine (2.1) and (2.2) to conclude
\[
^{u+t+1}a\, \equiv \, ^{y+t+1} a \mod L(t+1,n).
\]
The shown result is stronger than the claimed property, which may be obtained by putting $t=h$, $k=u+h$ and $j=y+h$.
\end{proof}

\section{Reducing $MTP$ to $IFP$}
In this section, we will discuss a polynomial-time reduction from tetration modulo natural numbers to computing the prime factorization of integers. We will use the results from the previous section and will follow the approach demonstrated in \cite[pp. 578-580]{BlaBor}.

Let $N\in\N$ and assume that we want to compute $^k a \pmod N$ for given $a,k\in\N$. For brevity, let $h=H(N)$. We start by computing the iterates of $\lambda(N)$. Assuming that we may compute the prime factorization of natural numbers in polynomial-time, this can be done by applying the formulae given in Lemma \ref{l1} iteratively on $N$, $\lambda(N)$, $\lambda^{(2)}(N)$ and so on, until we reach $\lambda^{(h)}(N)=1$. In (1) of Lemma \ref{l2}, we saw that $h\leq \lceil \log_2(N)\rceil$, hence the overall complexity for performing this procedure is polynomial in $N$. In the following, we assume knowledge of $h$ and the prime factorization of $N,\lambda(N)$, $\lambda^{(2)}(N),...,\lambda^{(h)}(N)$. To compute $^k a \pmod N$, we discuss two cases.
\vspace{12pt}

Case 1: $k>h$.

Theorem \ref{t1} implies that $^k a\, \equiv \, ^{h+1} a \mod N$ and, therefore, we have to find the smallest nonnegative integer $S$ such that $S\,\equiv\, ^{h+1} a \mod N$. We start by defining $B_1:=0$ and
\[
B_2:=\begin{cases} 1 &\mbox{if $a$ is odd},\\
0 & \mbox{if $a$ is even.} \end{cases}
\]
It is easy to see that $B_1\,\equiv \, ^1 a \mod \lambda^{(h)}(N)$ and $B_2\,\equiv \, ^2 a \mod \lambda^{(h-1)}(N)$. We now assume that $t\geq 2$ and that we know the unique integer $0\leq B_t <\lambda^{(h+1-t)}(N)$ satisfying
\[
B_t\,\equiv\, ^t a \mod \lambda^{(h+1-t)}(N).
\]
To derive $B_{t+1}$, we form the orthogonal decomposition 
\[
\lambda^{(h-t)}(N)=VW=V(a,\lambda^{(h-t)}(N))\cdot W(a,\lambda^{(h-t)}(N))
\]
of $\lambda^{(h-t)}(n)$ with respect to $a$. By definition, $\gcd(V,W)=\gcd(V,a)=1$ and every prime factor of $W$ also divides $a$. Therefore, one observes $^{t+1}a\, \equiv \, a^{B_t} \mod V$. Also note that
\[
^{t} a\geq 2^{t}\geq 2t=2H(\lambda^{(h-t)}(N))>E(\lambda^{(h-t)}(N))\geq E(W)
\]
yields $^{t+1} a\equiv 0 \mod W$. $B_{t+1}$ may now easily be found by an application of the Chinese Remainder Theorem. We proceed in this manner until we have computed the value $B_{h+1}$. Since $N=\lambda^{(0)}(N)=\lambda^{(h+1-(h+1))}(N)$, we may conclude that $S\equiv B_{h+1}\mod N$.
\vspace{12pt}

Case 2: $k\leq h$.

We first need to find the smallest number $g$ such that $^g a > N$. If $g$ is larger or equal to $k$, then the value  for $^k a\,\equiv\, a^{^{k-1}a} \mod N$ may be computed by applying the Square-and-Multiply procedure. Assume now that $k>g$. We first compute 
\[
B_g\, :=\,^g a\mod \lambda^{(k-g)}(N),
\]
also by using Square-and-Multiply. We then proceed inductively just like explained above. Note that we obtain $^{v-1} a>N$ for every $v>g$. We conclude that, for every $t=g+1,...,k$, it holds that
$^{v-1}a>E(W(a,\lambda^{(k-t)}(N)))$ and that
\[
^v a\,\equiv\, 0 \mod W(a, \lambda^{(k-t)}(N)).
\] 
This allows us to derive $B_t$ for $t=g+1,...,k$ by applying the Chinese Remainder Theorem via
\[
B_t\equiv\begin{cases} a^{B_{t-1}} &\mbox{$\mod V(a,\lambda^{(k-t)}(N))$},\\
0 & \mbox{$\mod W(a,\lambda^{(k-t)}(N))$.} \end{cases}
\]
Clearly, $B_k$ is the desired value. Note that, in both cases, (1) of Lemma \ref{l2} implies that the complete procedure can be performed in polynomial-time. 

\section{Reducing $IFP$ to $MTP$}
This section focuses on the main topic of this paper, namely the description of an efficient integer factorization method based on the assumption that a polynomial-time algorithm for $MTP$ is known.

We first roughly explain our core idea for the desired algorithm by taking a closer look on the terms $^{k+1}a\, -\, ^k a$. Indeed, we have
\begin{align*}
^{k+1}a\, -\, ^k a\, =\, ^ka(a^{^k a\, - \, ^{k-1} a}-1) \,=&\, ^ka(a^{^{k-1} a(a^{^{k-1}a\,-\,^{k-2} a}-1)}-1)\\
=& ^ka(a^{^{k-1} a(a^{^{k-2} a(a^{^{k-2}a\,-\,^{k-3}a}-1)}-1)}-1) 
\end{align*}
and so on. On each level of this tower-like expression, we have a term of the form $a^{^{k-i} a\, - \, ^{k-i-1} a}-1$. Now let $N\in\N$. Motivated by Lemma \ref{l3}, we consider these levels modulo the respective iterates of Carmichael's function $\lambda(N)$. If one of these strictly decreasing iterates divides the corresponding $^{k-i} a\, - \, ^{k-i-1} a$, this expression takes on the value $0$. In this case, a chain reaction is induced on the lower levels of the tower; they also become $0$, and the tower ``collapses" to a value divisible by $N$. If $N$ is composite, than there might be a $k$ for which this chain reaction happens modulo a prime divisor $p$ of $N$, but not modulo $N$. If this is the case, then $\gcd(^{k+1}a\, -\, ^k a,N)$ yields a nontrivial divisor.

We now discuss the details of our idea and start by considering the subsequent statement about the sequence of the values of tetration modulo natural numbers.

{\lemma{Let $n,a\in\N$. The sequence $(^k a \mod n)_{k\in\N_0}$ is periodic with period length $1$. Its preperiod length is smaller or equal to $H(n)+1$.}}

\begin{proof}
Since $n$ divides $L(n)$, this claim is an immediate consequence of Theorem \ref{t1}.
\end{proof}

\definition{Let $n,a\in\N$. By $\lev_n(a)$, the \emph{level of $a$ modulo $n$}, we denote the length of the preperiod of $(^k a \mod n)_{k\in\N_0}$.}
\vspace{12pt}

Suppose that we want to compute the prime factorization of $N\in\N$ under the assumption that we are able to compute modular tetration in polynomial-time. We choose some natural number $a\geq 2$. For $k=0,...,\lceil \log_2 N\rceil$, we then evaluate
\[
g_k:=\gcd(^{k+1} a\,-\, ^{k} a,N).
\]
Assume that $N$ is composite and there exists a prime factor $p$ of $N$ such that $t:=\lev_p(a)\neq \lev_N(a)$. By definition, this implies
\[
^{j+1} a\,\equiv \, ^{j} a \mod p
\]
for all $j\geq t$. However, there has to exist at least one $j\in\{t,...,\lceil \log_2 N\rceil\}$ such that the congruence above does \emph{not} hold modulo $N$; otherwise it would contradict the fact that $t<\lev_N(a)$. We conclude that, in this case, one of the values $g_k$ yields a nontrivial divisor of $N$. 

Obviously, the prime factorization of $N$ can be found in polynomial-time by reapplying the discussed algorithm to the found divisors. However, we made a crucial assumption for proving correctness of our method. Hence, we are left with solving the following computational problem. 

{\problem{Let $N\in\N$ be composite. Find any $a\geq 2$ such that there exists a prime factor $p$ of $N$ with $\lev_p(a)\neq \lev_N(a)$.}\label{p1}}
\vspace{12pt}

If we are able to construct a (deterministic) polynomial-time algorithm for this task, we clearly derive a (deterministic) polynomial-time reduction from $IFP$ to $MTP$. Therefore, the remainder of this paper focuses on its discussion. 

\example{Before we proceed, let us consider a few examples of applying the algorithm.
\begin{itemize}
\item{Let $N=60507095029$. We may compute $^6 2=2^{2^{2^ {16}}} \pmod N$ by squaring $2^{16}=65536$ times. In fact, one derives that $\gcd(^6 2\,-\,^5 2,N)=224951$ yields a nontrivial divisor of $N$.

Note that $k=6$ is the largest number such that $^k 2 \pmod N$ can be computed this way. For $k=7$, one would have to perform $2^{65536}$ squarings, which is infeasible.}
\item{
Consider the large semiprime number
\begin{align*}
N=&15226050279225333605356183781326374297180681149613\\
&80688657908494580122963258952897654000350692006139,
\end{align*} 
which has been labeled RSA-100 in the RSA-Factoring Challenge. Using the already known prime factorization
\begin{align*}
N&=37975227936943673922808872755445627854565536638199\\
 &\times 40094690950920881030683735292761468389214899724061
\end{align*}
of this number and the algorithm described in Section 3, one can show that $\gcd(^{10} 3\,-\,^9 3,N)$ is nontrivial. Considering all numbers up to $50$, the complete list of those bases $a$ for which we do \emph{not} find a factor of $n$ contains the $16$ elements $a=2,6,9,15,16,17,18,22,26,38,39,41,42,44,48,49$.
}
\item{The list of the bases up to $50$ for which we do not find a factor of RSA-110 is empty.}
\item{The list of the bases up to $50$ for which we do not find a factor of RSA-129 consists of $18$ elements.\label{e1}}
\end{itemize}}

\section{The level and iterated orders}
We continue by discussing results concerning bounds and expressions for $\lev_n(a)$ and will demonstrate its connection to the following notion of iterations of the order of $a$.

\definition{Let $n,a\in\N$. We define
\begin{itemize}
\item{$\ord_n^{(0)}(a):=n$,}
\item{$\ord_n^{(k)}(a):=\ord_{V(a,\ord_n^{(k-1)}(a))}(a)$ for all $k\in\N$,}
\end{itemize}
where $V$ is the respective factor of the orthogonal decomposition as it is described in Definition \ref{d1}. Furthermore, by $L_a(n)$ we denote the least common multiple of all the values $\ord_n^{(i)}(a)$ for $i\geq 0$.}

{\lemma{Let $n,a\in\N$. Then $\ord_n^{(k)}(a)\mid \lambda^{(k)}(n)$ for every $k\in\N_0$ and, as a result, $L_a(n)\mid L(n)$.
}\label{l9}}

\begin{proof}
We prove the statement by induction on $k$. The case $k=0$ is trivial. Now assume that the claim is true for $k-1$. We clearly have 
\[
a^{\lambda(V(a,\lambda^{k-1}(n)))}\equiv 1 \mod V(a,\lambda^{k-1}(n)).
\]
It follows from $V(a,\lambda^{k-1}(n))\mid\lambda^{k-1}(n)$ and Remark \ref{r1} that $\lambda(V(a,\lambda^{k-1}(n)))$ is a divisor of $\lambda(\lambda^{k-1}(n))=\lambda^{k}(n)$, and we derive
\[
a^{\lambda^{k}(n)}\equiv 1 \mod V(a,\lambda^{k-1}(n)).
\]
By the induction assumption, $\ord^{k-1}(a)$ divides $\lambda^{k-1}(n)$. It is easy to see that this yields 
\[
V(a,\ord^{k-1}(a))\mid V(a,\lambda^{k-1}(n)),
\]
hence  we conclude that the last congruence also holds modulo $V(a,\ord^{k-1}(a))$, which implies the claim.
\end{proof}

{\lemma{Let $n,a,k\in\N$. Then the following hold:
\begin{enumerate}
\item{If $^{k} a\, \equiv \, ^{k-1} a \mod n$, then $\ord_n^{(k)}(a)=1$,}
\item{$\lev_n(a)\geq \min\{\nu: \ord_n^{(\nu)}(a)=1\}-1$.}
\end{enumerate}\label{l4}}}

\begin{proof}
We first show (1). Assume that $^{k} a\, \equiv \, ^{k-1} a \mod n$. We have
\[
^{k-1}a(a^{^{k-1} a\, - \, ^{k-2} a}-1)\equiv 0 \mod n,
\]
and since $V(a,n)$ divides $n$, this implies $^{k-1} a\, \equiv \, ^{k-2} a \mod \ord_n^{(1)}(a)$. Next, we have
\[
^{k-2}a(a^{^{k-2} a\, - \, ^{k-3} a}-1)\equiv 0 \mod \ord_n^{(1)}(a)
\]
and easily obtain $^{k-2} a\, \equiv \, ^{k-3} a \mod \ord_n^{(2)}(a)$. We proceed in this manner until we reach
\[
a\,\equiv\,^1 a\,\equiv\, ^0 a\,\equiv\, 1 \mod\ord_n^{(k-1)}(a),
\]
and the statement follows.

For proving (2), let $l:=\lev_n(a)$. We derive $^{l+1} a\, \equiv \, ^{l} a \mod n$. By applying (1), it follows that $\ord_n^{(l+1)}(a)=1$. We therefore deduce 
\[
l+1\geq\min\{k: \ord_n^{(k)}(a)=1\},
\]
which we wanted to show.
\end{proof}

{\theorem{Let $n,a,k\in\N$ and suppose that $a$ is coprime to $L(n)$. Then $^{k} a\, \equiv \, ^{k-1} a \mod n$ holds if and only if $\ord_n^{(k)}(a)=1$.}}

\begin{proof}
Due to (1) of Lemma \ref{l4}, we only have to show that from $\ord_n^{(k)}(a)=1$ it follows that $^{k} a\, \equiv \, ^{k-1} a \mod n$. By definition, we have $a\equiv 1 \mod V(a, \ord^{(k-1)}(a))$. Since we assumed $\gcd(L(n),a)=1$, we also derive that $\gcd(a,\ord^{(j)}(a))=1$ for every $j=0,1,...,k$ and, hence,
$
a\equiv 1 \mod \ord^{(k-1)}(a).
$
Clearly, this implies $a^{a-1}-1\equiv 0 \mod V(a,\ord^{(k-2)}(a))$. By similar arguments, we may deduce 
\[
a^a(a^{a-1}-1)\equiv 0 \mod\ord^{(k-2)}(a).
\]
Considering the proof of (1) of Lemma \ref{l4}, it is easy to see that we may obtain $^{k} a\, \equiv \, ^{k-1} a \mod n$ by proceeding in this manner.
\end{proof}

As an immediate consequence of the results discussed above, we derive the following corollaries.

{\corollary{Let $n,a\in\N$ and suppose that $a$ is coprime to $L(n)$. Then 
\begin{align*}
\lev_n(a)&= \min\{\nu: \ord_n^{(\nu)}(a)=1\}-1\\
&=\min\{\nu:\ord_n^{(\nu)}(a)\mid a-1\}.
\end{align*}}\label{c4}}

{\corollary{Let $u,v\in\N$ and let $a\in\N$ be coprime to $L(\lcm(u,v))$. Then we have 
\[
\lev_u(a)= \lev_v(a)
\]
if and only if
\[
\forall j\geq 1: \lev_{\ord_u^{(j)}(a)}(a)= \lev_{\ord_v^{(j)}(a)}(a).
\]
}\label{c1}}

{\lemma{Let $n,a\in\N$ and $n=p_1^{e_1}\cdots p_l^{e_l}$ be the prime factorization of $n$. We have
\begin{enumerate}
\item{$\lev_n(a)=\max_{i\in\{1,...,l\}}\{\lev_{p_i^{e_i}}(a)\}$,}
\item{If $a$ is coprime to $L(n)$, then $\lev_n(a)=\lev_n(a\mod L_a(n))$.}
\end{enumerate}
}\label{l11}}

\begin{proof}
We start by showing $(1)$: Let $m:=\max_{i\in\{1,...,l\}}\{\lev_{p_i^{e_i}}(a)\}$. By definition, we obtain $\lev_n(a)\geq m$. Moreover, we have
$
^{j+1} a\,\equiv\, ^j a \mod p_i^{e_i}
$
for all $j\geq m+1$ and all $i\in\{1,...,l\}$. Since this implies that the congruence above also holds modulo $n$, we derive $\lev_n(a)\leq m$ and the claim follows.

For proving $(2)$, let $r:=a\pmod{L_a(n)}$ and note that $a\equiv r\mod L_a(n)$ clearly implies $\ord_n^{(i)}(a)=\ord_n^{(i)}(r)$ for all $i\geq 0$. Due to Corollary \ref{c4}, this already yields the desired result.
\end{proof}

{\lemma{Let $n,a\in\N$. Then $\lev_n(a)=\min\{k\in\N_0:\, ^{k+1}a\equiv \,^k a \mod n\}$.}\label{l13}}

\begin{proof}
Via induction on $k$, we prove that $^k a-\, ^{k-1}a$ divides $^{k+1} a-\, ^k a$ for every $k\in\N$. The result then follows by definition. For $k=1$, the claim is true. Assume that $^{k-1}a-\, ^{k-2} a\mid\, ^k a-\, ^{k-1}a$. We derive that $^k a-\, ^{k-1}a=\,^{k-1}a(a^{^{k-1}a-\,^{k-2}a}-1)$ divides $^{k}a(a^{^{k}a-\,^{k-1}a}-1)=\,^{k+1} a-\, ^{k} a$, which we wanted to show.
\end{proof}

We conclude this section by proving the following property of the level of elements modulo prime powers.

{\theorem{Let $p\in\PP$ and $a\in\N$ be coprime. To shorten notation, we define $M:=\min\{j\in\N: p\mid \ord_{p^j}(a)\}$. For all $n\in\N$, we then have
\[
\lev_{p^n}(a)=\lev_p(a)+\left\lceil \frac{n}{M-1}-1\right\rceil.
\]
}\label{t5}}

\begin{proof}
Let $n\in\N$ be arbitrary. We start with showing that
\[
\ord_{p^n}(a)=\begin{cases} \ord_p(a) &\mbox{for }  n<M,\\
p^{n-M+1}\cdot\ord_p(a) &\mbox{for } n\geq M. \end{cases}
\]
Let $b$ be a primitive root modulo $p^n$. Then there exists $r\in\Z$ coprime to $p$ such that $b^{rp^e}\equiv a \mod p^n$ for some $e\geq 0$. Furthermore, we have
\[
\ord_{p^n}(a)=\frac{p^{n-1}(p-1)}{\gcd(rp^e,p^{n-1}(p-1))}=\frac{p^{n-1}}{\gcd(p^e,p^{n-1})}\cdot \frac{p-1}{\gcd(r,p-1)}.
\]
We now show that $e=\min\{n-1,M-2\}$. If $n<M$, then $\ord_{p^n}(a)\mid p-1$ immediately implies $e=n-1$. For $n\geq M$, we have to prove that $e=M-2$. First assume that $e>M-2$. From $b^{rp^e}\equiv a \mod p^{M}$ and the fact that $b$ is also a primitive root modulo $p^M$, it follows that $p\nmid \ord_{p^M}(a)$, which contradicts the definition of $M$. Assuming $e<M-2$, we derive a contradiction to the minimality of $M$ by considering $b^{rp^e}\equiv a \mod p^{M-1}$. We conclude that $e=\min\{n-1,M-2\}$, and since $b^{p^e}$ is a primitive root modulo $p$ and $b^{rp^e}\equiv a \mod p$ holds, it also directly follows that $\ord_p(a)=(p-1)/\gcd(r,p-1)$. This implies the statement

Let $n\in\N$ and $j\in\N$ be arbitrary. We now prove the equality stated in the lemma. Using the assumption that $p$ and $a$ are coprime, one easily derives that $^{j+1}a\equiv\, ^j a\mod p^n$ is equivalent to $^ja\equiv\, ^{j-1} a\mod \ord_{p^n}(a)$. For $n<M$, it now follows from Lemma \ref{l13} that $\lev_{p^n}(a)$ is equal to the minimal value $m$ such that $^{m}a\equiv\, ^{m-1}a\mod \ord_p(a)$. It is clear that $m=\lev_p(a)$, which proves the claim. Now let $n=M$. According to the property shown above, we are looking for the minimal $m$ such that $^{m}a\equiv\, ^{m-1}a\mod p\cdot\ord_p(a)$. It is easy to check that $m=\lev_p(a)+1$. Indeed, we have seen that
\[
^{\lev_p(a)+1}a\equiv\, ^{\lev_p(a)}a\mod p^{M-1}\cdot\ord_p(a).
\]
As a consequence, $\lev_{p^n}(a)=\lev_p(a)+1$ for all $n\geq M$ satisfying $n-M+1\leq M-1$, hence for all $n\leq 2(M-1)$. Considering this argument for general $n\in\N$, we may conclude that $\lev_{p^n}(a)=\lev_p(a)+k$ for all $k(M-1)<n\leq (k+1)(M-1)$. This implies the equality.
\end{proof}

\section{Reducing $SDP$ to $MTP$}
We will now use the results from the last sections to prove that the problem of computing the squarefree part of a natural number is deterministically polynomial-time reducible to modular tetration.

{\definition{Let $N\in\N$. Then the \emph{squarefree part of $N$} is defined as the smallest number $r$ such that $N/r$ is a square.}}

\vspace{12pt}
The squarefree decomposition problem ($SDP$) is widely believed to be as hard as integer factorization itself. We will now first explain how to factorize non-squarefree numbers via modular tetration.

{\lemma{The problem of computing a nontrivial factor of a non-squarefree number is deterministically polynomial-time reducible to modular tetration.}}

\begin{proof}
Let $N\in\N$ be any non-squarefree number not divisible by $2$ or $3$ and $p$ a prime factor of $N$ such that $p^2\mid N$. As a consequence of Theorem 3 in \cite{Gra}, there is some $a\leq \log^2 p\leq \log^2 N$ such that $a$ is a $p$-th power nonresidue modulo $p^2$. For such $a$, we also derive that $p$ divides $\ord_{p^2}(a)$. Applying Theorem \ref{t5} with $n=M=2$ then yields that $\lev_{p^2}(a)=\lev_p(a)+1$. We conclude that there is some $a\leq \log^2 N$ such that $\lev_p(a)<\lev_N(a)$. Using the algorithm described in Section 4.4 for all $a\leq \log^2 N$, we will therefore eventually find a nontrivial factor of $N$.
\end{proof}

In the following proof, we will see how we may reduce $SDP$ to the nontrivial factorization of non-squarefree numbers and, hence, show our desired result.

{\theorem{$SDP$ is deterministically polynomial-time reducible to $MTP$.}\label{t3}}

\begin{proof}
In the preceding lemma, we have shown that a polynomial-time algorithm for $MTP$ allows us to efficiently find a nontrivial factor of any given non-squarefree number. Let $N$ now be the number of which we want to compute the squarefree part $r(N)$. We apply the discussed algorithm to $N$. If it fails to find a factor, $N$ itself is squarefree and $r(N)=N$. Assume now that we have found nontrivial $u,v$ such that $N=uv$. Note that for $g=\gcd(u,v)$, we clearly have $g^2\mid N$. Consequently, it is easy to see that
\[
r(N)=r(u/g)\cdot r(N/(ug)).
\]
Hence, the problem to determine $r(N)$ can be reduced to the problem of determining $r(u/g)$ and $r(N/(ug))$. We next apply our algorithm to the numbers $u/g$ and $N/(ug)$ to either determine that they are squarefree or to find nontrivial divisors. Imagine this procedure as the computation of a binary tree, where all nodes are divisors of $N$. We do not obtain further children of a node if the corresponding divisor is squarefree. Since the children nodes are always at least two times smaller than their parent nodes, the depth of this tree is bounded by $O(\log N)$. Furthermore, the product of all nodes on the same level of the tree is always smaller or at most equal to $N$, which implies the width of the tree is also bounded by $O(\log N)$. We conclude that the complete tree and the values of $r(\bullet)$ may be computed in polynomial-time.
\end{proof}

The importance of the squarefree decomposition problem is also due to the fact that one of the main task of computational algebra, namely determining the ring of integers of an algebraic number field, reduces to it deterministically and in polynomial-time. We derive the following corollary.

{\corollary{Determining the ring of integers of an algebraic number field is deterministically polynomial-time reducible to modular tetration.}\label{c5}}

\begin{proof}
This is a direct consequence of Theorem \ref{t3} and Theorem 4.4 in \cite{Len}.
\end{proof}

The technique discussed above demonstrates how to factor any natural number in squarefree divisors. It does not provide a rigorous method for the factorization of products of distinct primes. However, our calculations in Example \ref{e1} indicate the efficiency of the reduction from factoring large semiprimes to modular tetration. In the following section, we will provide a heuristic analysis of this reduction.

\section{Analysis of reducing $IFP$ to $MTP$}
Let $N\in\N$ be the number which we want to factorize completely. We have already seen in the last section that we may restrict our attention to the case where $N$ is squarefree and divisible by at least two distinct prime factors $p$ and $q$. Problem \ref{p1} then reduces to finding $a\in\N$ such that $\lev_p(a)\neq\lev_q(a)$. We will now continue to investigate this situation and start by introducing the following notation.

\definition{Let $u,v,a\in\N$. We define $k_{u,v}:=\lcm(u,v)$ and the set
\[
W_{u,v}:=\{a\in\Z_{L(k_{u,v})}^*: \lev_u(a)=\lev_v(a)\}.
\]
Moreover, let $\omega(u,v):=|W_{u,v}|/\phi(L(k_{u,v}))$, where $\phi$ denotes the totient function.}
\vspace{12pt}

Our goal is to find an upper bound for $\omega(p,q)$ for distinct primes $p$ and $q$. The following lemma will be useful in the subsequent proofs.

{\lemma{Let $a,b,c\in\N$ such that $a\mid b$ and $\gcd(a,c)=1$. Then
\[
|\{x\in\Z_b^*: x\equiv c \mod a\}|=\frac{\phi(b)}{\phi(a)}.
\]
}\label{l5}}

\begin{proof}
Let $b=W(a,b)\cdot V(a,b)$ be the orthogonal decomposition of $b$. Moreover, write $b=p_1^{\beta_1}\cdots p_l^{\beta_l}$ and $a=\prod_{i\in M} p_i^{\alpha_i}$ for $\alpha_i\leq \beta_i$ and some $M\subseteq\{1,...,l\}$. We have $W(a,b)=\prod_{i\in M}p_i^{\beta_i}$.

Let $p_\nu$ be an arbitrary prime factor of $W(a,b)$. We show that there are exactly $p_\nu^{\beta_\nu-\alpha_\nu}$ elements $x\in\Z_{p_\nu^{\beta_\nu}}^*$ such that $x\equiv c\mod p_\nu^{\alpha_\nu}$. Clearly, there is one such element in $\Z_{p_\nu^{\alpha_\nu}}^*$, namely $c$ itself. Modulo $p_\nu^{\beta_\nu}$, we then have the elements
\[
p_\nu^{\alpha_\nu}+c, 2p_\nu^{\alpha_\nu}+c,\ldots,(p_{\nu}^{\beta_\nu-\alpha_\nu})p_\nu^{\alpha_\nu}+c.
\]
The Chinese Remainder Theorem now yields that there are $\prod_{i\in M} p_i ^{\beta_i-\alpha_i}$ elements $x\in\Z_{W(a,b)}^*$ such that $x\equiv c \mod a$. By one further application of CRT on those elements and the members of $\Z_{V(a,b)}^*$, we may derive the desired result by noting that
\[
\phi(V(a,b))\cdot \prod_{i\in M} p_i ^{\beta_i-\alpha_i}=\phi(b)/\phi(a).
\]
\end{proof}

{\lemma{Let $u,v\in\N$ and set $k=k_{u,v}$. If $k$ is coprime to $L(\lambda(k))$, then we have
\[
\omega(u,v)= \frac{1}{\phi(k)}\cdot\sum_{a\in \Z^*_{k}} \omega(\ord_u(a),\ord_v(a)).
\]
}\label{l10}}

\begin{proof}
Let $a\in\Z^*_{L(k)}$ be arbitrary. As a consequence of Corollary \ref{c1}, we derive that $\lev_u(a)=\lev_v(a)$ if and only if $\lev_{\ord_u(a)}(a)=\lev_{\ord_v(a)}(a)$. Therefore, $a\in W_{u,v}$ if and only if $a\pmod {L(\lambda(k))}\in W_{\ord_u(a),\ord_v(a)}$. 

Now fix any $\alpha\in\Z^*_{k}$ and set $o_u=\ord_u(\alpha)$ and $o_v=\ord_v(\alpha)$. By definition, the members of $W_{o_u,o_v}$ are also elements in $\Z^*_{L(k_{o_u,o_v})}$. We consider the cardinality of the set
\[
C_\alpha=\{x\in\Z_{L(\lambda(k))}^*: \lev_{o_u}(x)= \lev_{o_v}(x)\}.
\]
Since $L(k_{o_u,o_v})\mid L(\lambda(k))$ holds, we may apply Lemma \ref{l5} to derive that $|C_\alpha|$ is equal to
\[
|W_{o_u,o_v}|\cdot\frac{\phi(L(\lambda(k)))}{\phi(L(k_{o_u,o_v}))}=\phi(L(\lambda(k)))\cdot \omega(o_u,o_v).
\]
Note that $k$ and $L(\lambda(k))$ are coprime. Due to the Chinese Remainder Theorem, for each $x\in C_\alpha$ there is exactly one element $r\in\Z_{L(k)}^*$ satisfying
\begin{align*}
r&\equiv \alpha \mod k,\\
r&\equiv x \mod L(\lambda(k)).
\end{align*}
As a result of these observations, we derive
\[
|W_{u,v}|= \sum_{a\in \Z^*_{k}} |C_{a}|= \phi(L(\lambda(k)))\cdot\sum_{a\in \Z^*_{k}} \omega(\ord_u(a),\ord_v(a)).
\]
We conclude that
\[
\omega(u,v)=\frac{|W_{u,v}|}{\phi(L(k))}=\frac{1}{\phi(k)}\cdot\sum_{a\in \Z^*_{k}} \omega(\ord_u(a),\ord_v(a)),
\]
where we have used $\gcd(k,L(\lambda(k)))=1$ again.
\end{proof}

For primes $p$ and $q$, we are now ready to prove the following upper bound for the value of $\omega(p,q)$.

{\theorem{Let $p,q$ be primes such that $p< q$ and set $n:=pq$. Then
\[
\omega(p,q)\leq\frac{1}{\phi(n)}\cdot\sum_{a\in \Z^*_{n}} \omega(\ord_p(a),\ord_q(a)),
\]
where equality holds if $p\nmid L(q)$.
}\label{t4}}

\begin{proof}
Let $a\in\Z_{L(n)}^*$ be arbitrary. Assume that there is $j\geq 1$ with $p\mid \ord_q^{(j)}(a)$, then we derive $\lev_q(a)>\lev_{\ord_q^{(j)}(a)}(a)$ by Corollary \ref{c4}. Moreover, we have $\lev_{\ord_q^{(j)}(a)}(a)\geq \lev_p(a)$ by definition. As a result, we deduce $\lev_q(a)\neq \lev_p(a)$ and $a\not\in W_{p,q}$. Consequently, $W_{p,q}$ is a subset of $\Z_{L(n)}^*\backslash M$, where
\[
M:=\{a\in\Z_{L(n)}^*\mid\exists j\geq 1: p\mid \ord_q^{(j)}(a)\}.
\]
It is easy to see that $M$ is empty if and only if $p\nmid L(q)$. In this case, we have $\gcd(n,L(\lambda(n)))=1$ and we may apply Lemma \ref{l10} directly to deduce the claim. Suppose that $M\neq \emptyset$ and $\gcd(n,L(\lambda(n)))=p$. We follow the proof of Lemma \ref{l10}. Fix any $\alpha\in\Z^*_{n}$ and set $o_p=\ord_p(\alpha)$ and $o_q=\ord_q(\alpha)$. For any element $x$ in the set $C_\alpha$ defined in the proof of Lemma \ref{l10}, let $\mathcal{C}(\alpha,x)$ denote the element $r\in\Z_{L(n)}^*$ for which
\begin{align*}
r&\equiv \alpha \mod n,\\
r&\equiv x \mod L(\lambda(n)).
\end{align*}
If such $r$ should not exist, we consider the label $\mathcal{C}(\alpha,x)$ as empty. Our goal is to estimate the cardinality of the set
\[
B_\alpha:=\{\mathcal{C}(\alpha,x): x\in C_\alpha\}.
\]
One easily checks that we have $B_\alpha\subseteq W_{p,q}\subseteq \Z_{L(n)}^*\backslash M$ and $W_{p,q}=\bigcup_{a\in \Z^*_{n}} B_a$.

We first define
\[
R_\alpha:=\{x\in\Z_{\kappa}^*: \lev_{o_p}(x)= \lev_{o_q}(x)\},
\]
where $\kappa:=L(\lambda(n))/p^e$ and $e:=\nu_p(L(\lambda(n)))\geq 1$ is the $p$-adic valuation of $L(\lambda(n))$. Moreover, for any $x\in Z_{L(\lambda(n))}^*$, let $x_\kappa:=x\pmod\kappa$. 

We now consider $\beta\in B_\alpha$. Since $x:=\beta\mod L(\lambda(n))\in C_\alpha$, by definition we derive $\lev_{o_p}(x)=\lev_{o_q}(x)$. Noting $\beta\not\in M$ and that $p$ does not divide $\ord_q^{(j)}(\beta)$ for any $j\geq 1$, Lemma \ref{l11} $(2)$ yields $\lev_{o_p}(x)=\lev_{o_p}(x_\kappa)$ and $\lev_{o_q}(x)=\lev_{o_q}(x_\kappa)$. Obviously, this implies that $\beta\pmod \kappa\in R_\alpha$. As a consequence, we obtain that $\beta=\mathcal{C}(\alpha,x)$, where $x$ satisfies $x_\kappa \in R_\alpha$. Due to the fact that $\beta\equiv x\mod L(\lambda(n))$, we have $\ord_q^{(j)}(\beta)=\ord_{o_q}^{(j-1)}(\beta)=\ord_{o_q}^{(j-1)}(x)$ for all $j\geq 1$. Therefore, $p\nmid \ord_{o_q}^{(j)}(x)$ for every $j\geq 0$. We deduce that $L_x(o_q)$ is a divisor of $\kappa$.

Since $\beta$ was arbitrary, these observations yield
\[
B_\alpha \subseteq \{\mathcal{C}(\alpha,x)\mid x\in C_\alpha '\},
\]
where $C_\alpha':=\{x\in C_\alpha: x_\kappa \in R_\alpha\wedge L_x(o_q)\mid \kappa\}$. We now want to prove that 
\[
C_\alpha'=\{x\in \Z_{L(\lambda(n))}^*: x_\kappa \in R_\alpha\wedge L_{x_\kappa}(o_q)\mid \kappa\}.
\]
Let $x\in C_\alpha'$, then $x_\kappa\in R_\alpha$ and $L_x(o_q)\mid \kappa$. We have to show that $L_{x_\kappa}(o_q)\mid \kappa$. In order to do this, we prove $L_x(o_q)=L_{x_\kappa}(o_q)$ by showing $\ord_{o_q}^{(j)}(x)=\ord_{o_q}^{(j)}(x_\kappa)$ via induction on $j$. For $j=0$, the statement is clear. Let the claim be true for $j-1$. Then
\begin{align*}
\ord_{o_q}^{(j)}(x)=\ord_{\ord_{o_q}^{(j-1)}(x)}(x)=&\ord_{\ord_{o_q}^{(j-1)}(x_\kappa)}(x)\\
=&\ord_{\ord_{o_q}^{(j-1)}(x_\kappa)}(x_\kappa)=\ord_{o_q}^{(j)}(x_\kappa),
\end{align*}
where we have used $\ord_{o_q}^{(j-1)}(x_\kappa)=\ord_{o_q}^{(j-1)}(x)\mid\kappa$ by applying the induction assumption and $L_x(o_q)\mid \kappa$.

Let $x$ be an element of the set on the right-hand side, then $x_\kappa\in R_\alpha$ and $L_{x_\kappa}(o_q)$ divides $\kappa$. We have to show that $x\in C_\alpha$ and $L_x(o_q)\mid \kappa$. For the latter one, we may again prove that $L_x(o_q)=L_{x_\kappa}(o_q)$ by using similar arguments as above. Now if we already know that $L_x(o_q)\mid \kappa$, we clearly have $\lev_{o_q}(x)=\lev_{o_q}(x_\kappa)$. And since it is easy to see that $L_x(o_p)$ divides $\kappa$, we also have $\lev_{o_p}(x)=\lev_{o_p}(x_\kappa)$. We finally deduce $x\in C_\alpha$, which concludes the proof of the claim that the sets are equal. Our next task is to prove that, for every residue $m\in\Z_{p}^*$, we have
\[
|\{x\in C_\alpha': x\equiv m \mod p\}|=\frac{|C_\alpha'|}{p-1}.
\]
For every $\xi\in R_\alpha$ satisfying $L_{\xi}(o_q)\mid \kappa$, there are $\phi(p^e)$ elements in $\Z_{p^e}^*$ which are of the form $\xi+j\kappa$ for $j\in\Z$. Clearly, the residues of these elements modulo $p$ are equally distributed among the numbers $1,2,...,p-1$. As a result, the same holds for the elements of the set $C_\alpha'$.

Following the concluding argument in the proof of Lemma \ref{l10} we observe that, due to the Chinese Remainder Theorem, $\mathcal{C}(\alpha,x)$ exists only for those $x\in C_\alpha'$ satisfying $x\equiv \alpha \mod p$. We know that there are $|C_\alpha'|/(p-1)$ such elements in $C_\alpha'$. We henceforth deduce
\begin{align*}
|W_{p,q}|=\sum_{a\in \Z^*_{n}} |B_a|\leq \sum_{a\in \Z^*_{n}} \frac{|C_a'|}{p-1}\leq& \sum_{a\in \Z^*_{n}} \frac{|C_a|}{p-1}\\
=&\frac{\phi(L(\lambda(n)))}{p-1}\cdot\sum_{a\in \Z^*_{n}} \omega(\ord_p(a),\ord_q(a)).
\end{align*}
Here, we have used the expression for the cardinality of $C_a$, which has been shown in the proof of Lemma \ref{l10} without using the assumption that $n$ and $L(\lambda(n))$ are coprime. Since it is easy to prove that $\gcd(n,L(\lambda(n)))=p$ implies 
\[
\frac{\phi(L(\lambda(n)))}{\phi(L(n))(p-1)}=\frac{1}{\phi(n)},
\]
we derive the desired bound.
\end{proof}

{\corollary{Let $p,q$ be primes such that $p<q$ and set $n:=pq$. Then
\[
\omega(p,q)\leq\frac{1}{\phi(n)}\cdot\sum_{r\mid p-1}\phi(r)\sum_{s\mid q-1} \phi(s)\cdot\omega(r,s),
\]
where equality holds if $p\nmid L(q)$.
}\label{c2}}
\begin{proof}
Note that, for any prime number $\rho$ and any divisor $d$ of $\rho-1$, there are exactly $\phi(d)$ elements in $Z_\rho^*$ of order $d$. Therefore, for any pair $(r,s)\in\Z_p^*\times\Z_q^*$ such that $r\mid p-1$ and $s\mid q-1$, the Chinese Remainder Theorem yields that we have exactly $\phi(r)\cdot\phi(s)$ elements $a\in \Z^*_{n}$ such that $\ord_p(a)=r$ and $\ord_q(a)=s$. This already proves the statement.
\end{proof}

We now discuss observations concerning the bound in Corollary \ref{c2}. Let $p$ and $q$ be distinct primes and assume $p<q$. First note that the bound is always smaller or equal to $1$ and, in the most usual case where $p\nmid L(q)$, is in fact an exact expression for $\omega(p,q)$. Two further points are worth to mention:
\begin{enumerate}
\item{One easily proves that $\omega(r,s)=1$ if and only if $r=s$. Considering our bound above, $\omega(p,q)$ tends to be larger if there are a lot of pairs $(r,s)$ such that $r\mid p-1$, $s\mid q-1$ and $r=s$. Now the number of these pairs obviously depends on $g=\gcd(p-1,q-1)$. For this reason, integers like $N=1541=23\cdot 67$ are much harder to factorize via modular tetration than average numbers. In this example, we have $g=22$ and $\omega(23,67)\approx 0.82$. Calculations have shown that $\omega$ can be estimated to be even closer to $1$ for larger examples with the same property.

However, we want to stress that there are other efficient methods to factor $N=pq$  in the exceptional case where $g$ is a large divisor of or even equal to $p-1$. Let $a\in Z_N^*$ and note that $g\mid N-1$. Under our assumption, $p-1$ divides $l(N-1)$ for $l$ being a product of some small primes or even for $l=1$. As a consequence, there is a high chance to obtain
\[
p=\gcd(a^{l(N-1)}-1,N).
\]}
\item{If $r,s$ are distinct primes, the argument in (1) can be applied to $\omega(r,s)$. More generally, one may observe that $\omega(p,q)$ tends to be smaller (larger) if
\[
\gcd(\lambda^{(k)}(p),\lambda^{(k)}(q))
\]
is relatively small (large) for the first few values of $k$.}
\end{enumerate}

\section{Further observations}
We conclude this paper by considering problems and open questions concerning the relation between certain algorithmic problems. We start by showing the following result.

{\theorem{Modular tetration is deterministically polynomial-time reducible to the discrete logarithm problem.}\label{t2}}

\begin{proof}
Let $N,a\in\N$ and assume that we want to compute $^k a\pmod N$ for any $k\in\N_0$. Set $h=\lceil \log_2 N\rceil$ and consider the algorithm described in Section 3. We will show that it suffices to know the iterates of the order of $a$ modulo $N$ instead of the iterates of Carmichael's function. Since a method for solving discrete logarithms allows to efficiently compute $\ord_N^{(i)}(a)$ for every $i\in\{0,...,h+1\}$, this will prove the statement.

First note that Theorem \ref{t1} and Lemma \ref{l2} imply that $^k a\equiv\, ^{h+1} a\mod N$. Hence, we have to compute $^m a\mod N$ for $m:=\max\{k,h+1\}$. In order to do so, we follow the approach of Case 2 in the algorithm of Section 3. We first compute the smallest number $g$ such that $^g a>N$. If $g\geq m$, we may compute our desired value $^m a\pmod N$ by using Square-and-Multiply. Assume that $m>g$. We then compute
\[
B_g:=\,^g a \mod \ord_N^{(m-g)}(a),
\]
also via Square-and-Multiply. We have $\ord_N^{(m-g)}(a)=\ord_{V(a,\ord_N^{(m-g-1)}(a))}(a)$
by definition. As a consequence of the fact that $a$ is coprime to $V(a,\ord_N^{(m-g-1)}(a))$, we immediately derive 
\[
^{g+1}a\equiv a^{B_g}\mod V(a,\ord_N^{(m-g-1)}(a)).
\]
Now since $^g a>N$ holds, we derive $^{g+1} a\equiv 0 \mod W(a, \lambda^{(m-g-1)}(N))$ by the same means as in Section 3. Due to Lemma \ref{l9}, $\ord_N^{(m-g-1)}(a)$ divides $\lambda^{(m-g-1)}(N)$ and, hence, we also have
\[
^{g+1} a\equiv 0 \mod W(a, \ord_N^{(m-g-1)}(a)).
\]
Note that we obtain $V$ and $W$ by computing $W=\gcd\left(a^{\lceil \log_2 N\rceil},\ord_N^{(m-g-1)}(a)\right)$. The value of $^{g+1} a\pmod{\ord_N^{(m-g-1)}(a)}$ may now be computed by an application of the Chinese Remainder Theorem. 

Repeating the computation discussed above, we conclude that we may compute the value of $^m a\pmod{N}$ in polynomial-time.
\end{proof}

\remark{The following graph represents the relations between $IFP$, $MTP$ and two other algorithmic problems discussed in this paper. Regarding the labels of the nodes, $CLP$ stands for the problem of computing Carmichael's function and $DLP$ for the discrete logarithm problem. A single-lined (double-lined) arrow from Problem A to Problem B indicates that Problem A is probabilistically (deterministically) polynomial-time reducible to Problem B. The dashed lines leading to $MTP$ indicate the efficient reduction from $IFP$ to $MTP$ and, consequently, from $CLP$ to $MTP$.

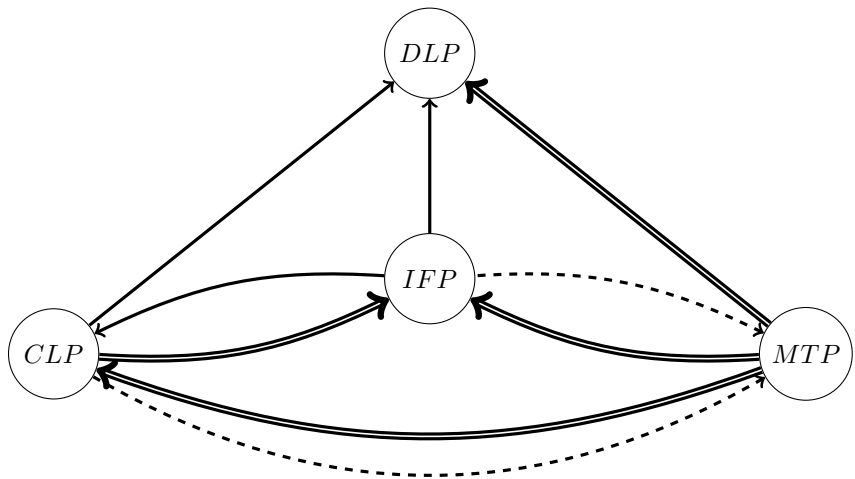
\begin{tikzpicture}
\node[minimum size=1.2cm] (A) at (10,0) [circle,draw] {$MTP$};
\node[minimum size=1.2cm] (B) at (5,1) [circle, draw] {$IFP$};
\node[minimum size=1.2cm] (C) at (0,0) [circle,draw] {$CLP$};
\node[minimum size=1.2cm] (D) at (5,4) [circle, draw] {$DLP$};

\draw[<-, dashed, very thick] (A) to[bend right=15]
   node[right] {} (B);
\draw[double, <-, very thick] (B) to[bend right=15]
   node[right] {}  (A);
\draw[<-, very thick] (D) to (C);
\draw[->, very thick] (B) to
   node[right] {} (D);
\draw[double, ->, very thick] (A) to
   node[right] {} (D);
\draw[double, ->, very thick] (C) to[bend right=15]
   node[right] {} (B);
\draw[->, very thick] (B) to[bend right=15]
   node[right] {} (C);
\draw[double, <-, very thick] (C) to[bend right=20]
   node[right] {} (A);
\draw[<-, dashed, very thick] (A) to[bend left=30]
   node[right] {} (C);
\end{tikzpicture}

Most of these relations are well-known and/or easy to check, and some of them have been shown or discussed in previous sections.}
\vspace{12pt}

We now state two problems which remained unsolved or undiscussed.

\problem{Show that $IFP$ is polynomial-time reducible to $MTP$.\label{p2}}
\vspace{6pt}

Due to Theorem \ref{t2}, the yet open problem of proving that there exists a deterministic polynomial-time reduction from $IFP$ to $DLP$ reduces to providing a deterministic solution for Problem \ref{p2}.

Of course, the results of this paper also show that there is a great interest in the following computational problem.

\problem{Find a suitable method for solving instances of modular tetration.}
\vspace{6pt}

If we are able to compute tetration modulo natural numbers of certain shape, for some bases $a$ or, more generally, for any specific set of input triples $(a,k,N)$, such method may give rise to special-purpose integer factorization algorithm.


\end{document}